\pdfoutput=1

\documentclass{article}
\usepackage[T1]{fontenc}
\usepackage[UKenglish]{babel}

\usepackage{Factory, Math, Theorema, Classic, Styl}

\DeclareMathAlphabet{\mathpzc}{OT1}{pzc}{m}{it}

\usepackage[all, 2cell]{xy} \UseAllTwocells 
\newcommand{\dotar}{\ar@{.>}}
\newcommand{\isoar}{\ar@{=} }			
\newcommand{\upar}{ \ar@<+.6ex> }
\newcommand{\car}{\ar@/^1pc/}

\newcommand{\B}{\mathbf{B}}
\newcommand{\op}{\text{op}} 
\newcommand{\stir}[2]{\genfrac{ \{ }{ \} }{0pt}{}{#1}{#2}}
\newcommand{\freenum}[2][\Z]{#1\binom{#2}{-}}

\newcommand{\de}{\diamond} 					
\DeclareMathOperator*{\De}{ \lozenge } 
\newcommand{\dev}[3]{#1(#2\de\cdots\de #3)}   	
\newcommand{\Dev}[2]{#1\left(\De #2\right)}

\newcommand{\Sets}{\mathfrak{Set}}
		
\newcommand{\Mod}{\mathfrak{Mod}} 
\newcommand{\FMod}{\mathfrak{FMod}} 		
\newcommand{\AFMod}{{}_A\FMod} 			
\newcommand{\AMod}{{}_A\Mod} 				
 			
\newcommand{\SMod}{{}_S\Mod} 	
		
\newcommand{\CAlg}{\mathfrak{CAlg}} 	

\newcommand{\Num}{\mathpzc{N}}							
\newcommand{\HPol}{\mathpzc{P}}

\newcommand{\BM}{P_n(M)}					
\newcommand{\BHK}[1][-]{P_n\Hom(K,#1)}	
\newcommand{\BH}[1][-]{P_n\Hom(\B^n,#1)}
\newcommand{\BB}{P_n(\B^{n\times n})}				
\newcommand{\BBMod}{{}_{\BB}\Mod}

\newcommand{\GB}{\Gamma^n(\B^{n\times n})} 
\newcommand{\GM}{\Gamma^n(M)} 
\newcommand{\GH}[1][-]{\Gamma^{n}\Hom(\B^n,#1)}
\newcommand{\GBMod}{{}_{\GB}\Mod}

\usepackage{fancyhdr}
\begin{document}

\lefthyphenmin=2 \righthyphenmin=2

\titul{POLYNOMIAL FUNCTORS \\ OF MODULES}
\auctor{Qimh Richey Xantcha%
\thanks{\textsc{Qimh Richey Xantcha}, Stockholm University: \texttt{qimh@math.su.se}}}
\datum{\today}
\maketitle

\bigskip 
\epigraph{\begin{vverse}[Och när jag stod där gripen, kall av skräck]
Och när jag stod där gripen, kall av skräck \\
och fylld av ängslan inför hennes tillstånd \\
begynte plötsligt mimans fonoglob \\
att tala till mig på den dialekt \\
ur högre avancerad tensorlära \\
som hon och jag till vardags brukar mest. \par 
\vattr{Harry Martinson, \emph{Aniara}}
\end{vverse}}

\bigskip 

\begin{argument} \noindent
We introduce the notion of numerical functors to generalise Eilenberg \& MacLane's 
polynomial functors to modules over a binomial base ring. 
After shewing how these functors are encoded by modules over a certain ring, 
we record a precise criterion
for a numerical (or polynomial) functor to admit a strict polynomial structure in the sense 
of Friedlander \& Suslin. 
We also provide several characterisations of analytic functors. 

\end{argument}

\bigskip 

\noindent 
It will be recalled that \emph{polynomial functors} were invented by Eilenberg \& MacLane 
\cite{EM} in 1954, and  
\emph{strict polynomial functors} by Friedlander \& Suslin \cite{FS} in 1997. 
They have since found numerous applications in algebraic topology and constitute to-day 
an active field of research. 

As evinced by terminology, the notion of polynomial functor is weaker than that 
of strict polynomial functor. It will naturally be enquired: how much weaker? 
The purpose of the present note is to provide a satisfactory answer to the question: 
\emph{When is a polynomial functor strict polynomial?} 

Consider then a commutative, unital ring $\B$, and 
let ${}_\B\Mod$ be the category of left $\B$-modules. 
Our first result shews that, when considering module functors on ${}_\B\Mod$, 
it will be no great restriction to consider only functors defined on 
${}_\B\FMod$, the subcategory of finitely generated and free modules:

\begin{inttheorem}[\ref{S: Bouc}]
Any functor $F\colon {}_\B\FMod\to{}_\B\Mod$ has a unique extension to a 
functor $\hat F\colon {}_\B\Mod\to{}_\B\Mod$ that is right-exact 
and commutes with filtered inductive limits.
\end{inttheorem}

Polynomial functors were initially conceived for abelian groups. While the notion, 
as such, is perfectly sensible for modules over any ring, and indeed admits a wide range of applications, 
scalar multiplication is nowhere taken into account. As a remedy for this, we introduce the notion 
of \emph{numerical functor} (Definition \ref{D: Numerical}), designed to make sense for any 
base ring $\B$ that is \emph{binomial}, which is to say equipped with binomial co-efficients. 
Since polynomial and numerical functors concur over $\Z$, all our results 
will be valid for integral polynomial functors. 

There is a characterisation of numerical functors that is found to tie in closely with strict polynomial functors:

\begin{inttheorem}[\ref{S: BAlg-polynomial}]
The functor $F\colon{}_\B\FMod\to{}_\B\Mod$ is numerical of degree $n$ 
if and only if its arrow maps extend to a system  
$$
F_A \colon A\otimes_\B \Hom_\B(M,N) \to A\otimes_\B \Hom_\B(F(M),F(N))
$$
of maps of degree $n$, multiplicative and natural in the binomial $\B$-algebra $A$.
\end{inttheorem}

We explore the elementary properties of numerical functors of degree $n$.
They constitute an abelian category $\Num_n$, which 
will be found Morita equivalent to the category of modules over 
a certain ring $\BB$, with 
$\B^{n\times n}$ betokening the $n\times n$ matrix ring of $\B$, 
and $P_n$ a construction given in Section \ref{A: Quasi-homogeneous}:

\begin{inttheorem}[\ref{S: Num}] 
The functor $P_n\Hom_\B(\B^n,-)$ is a small projective generator for $\Num_n$, through 
which there is an equivalence of categories
$$ 
\Phi\colon \Num_n \to \BBMod, \qquad F\mapsto F(\B^n).
$$ 
The left module structure on $F(\B^n)$ is given by $ [\tau]\cdot x = F(\tau)(x)$. 
\end{inttheorem}

The theorem was first enunciated for polynomial functors by Pirashvili \cite{Pirashvili}.
As a moral consequence, all theorems valid for 
integral polynomial functors will remain so, \emph{mutatis mutandis}, for numerical functors.

Suslin, Friedlander \& Bendel \cite{SFB} gave us the correspondence between homogeneous 
(strict polynomial) functors of degree $n$, which form an abelian category $\HPol_n$, and modules over the 
\emph{Schur algebra} $\GB$. This is Theorem \ref{S: HPol} below. 
A link is provided by the \emph{divided power map}
$$ 
\gamma_n\colon \BB\to\GB, \qquad [\sigma]\mapsto \sigma^{[n]}, 
$$
which, by restriction of scalars, gives rise to a forgetful 
functor from homogeneous to numerical functors: 

\begin{inttheorem}[\ref{T: Forgetful}]
The forgetful functor $\HPol_n\to\Num_n$ corresponds to restriction of scalars 
$$
\GBMod\to\BBMod
$$ 
along the divided power map $\gamma_n \colon \BB \to \GB$.
This functor is exact and faithful. It is full and essentially injective when $n=0,1,2$.
\end{inttheorem}

In Example \ref{X: gamma_3} (p.~\pageref{X: gamma_3}), we observe that the forgetful functor 
is not, in general, full or essentially injective for $n\geq 3$.

We then answer the question of when a numerical functor admits a strict polynomial structure. 
After introducing the concept of \emph{quasi-homogeneous} functors in Section \ref{A: Quasi-homogeneous} 
and shewing how such functors correspond to modules over $\Im\gamma_n$ (Theorem \ref{S: Quasi-homogeneous}), 
we prove:

\begin{inttheorem}[\ref{T: PFT}: The Polynomial Functor Theorem] 
Let $F$ be a numerical functor of degree $n$. Then $F$ may be given the structure of homogeneous 
functor of degree~$n$ if and only if $F$ is quasi-homogeneous of degree $n$ and the 
 $\Im\gamma_n$-module structure on $F(\B^n)$ may be extended to 
a $\GB$-module structure. 
\end{inttheorem}

We proceed to discuss analytic functors in some detail. An \emph{analytic} functor may be defined as 
a family of module functors $E_A\colon {}_A\FMod\to{}_A\Mod$, with $A$ ranging over all binomial $\B$-algebras, 
commuting with extension of scalars. Equivalently (Theorem \ref{S: ABC}), it is a functor
$J\colon {}_\B\Mod\to{}_\B\Mod$ with arrow maps  
$$
J_A\colon \Hom_A(A\otimes_\B M,A\otimes_\B N) \to \Hom_A(A\otimes_\B J(M),A\otimes_\B J(N)),
$$ 
multiplicative and natural in the binomial algebra $A$. We prove:

\begin{inttheorem}[\ref{S: Analytic}]
The analytic functors are precisely the filtered inductive limits of numerical functors. 
\end{inttheorem}


A combinatorial approach to our subject is discussed in the paper \cite{Laby}.

\section{Module Functors}

For the entirety of this article, $\B$ shall denote a fixed base ring of scalars, assumed 
commutative and unital.
All modules, algebras, homomorphisms and tensor products shall be taken over this 
$\B$, unless otherwise stated. 
We let $\Mod={}_\B\Mod$ denote the category of (unital) modules over this ring, and $\FMod$ the subcategory 
of finitely generated and free modules. 

A \textbf{module functor} shall be taken to mean a functor $\FMod \to \Mod$, 
usually non-linear.  
We shall be wholly content to consider such restricted functors exclusively, as is customarily done, 
for, as we presently expand upon, 
a functor defined on the subcategory $\FMod$ always has a canonical well-behaved extension 
to the whole module category

First, let us recall from Bouc \cite{Bouc} 
that a (non-additive) functor $F$ between abelian categories is \textbf{right-exact} 
if for any exact sequence 
$$ 
\xymatrix{ A\ar[r]^{\alpha} & B\ar[r]^{\beta} & C\ar[r] & 0, } 
$$
the associated sequence 
$$ 
\xymatrix{ F(A\oplus B)\ar[rr]^-{\substack{ F(\alpha+1_B) -F(1_B)}} 
& & F(B)\ar[r]^{F(\beta)} & F(C)\ar[r] & 0 }
$$
is also exact.  
This definition agrees with the usual one in the case of an additive functor. 
In fact, the usual definition actually \emph{implies} additivity of the functor, 
which renders it useless for our purposes. 

A variant of the following result, in the abstract setting of abelian categories, 
is due to Bouc (\cite{Bouc}, Theorem 2.14). 

\begin{theorem} 			\label{S: Bouc}
Any functor $F\colon \FMod\to\Mod$ has a unique extension to a 
functor $\hat F\colon \Mod\to\Mod$ that is right-exact 
and commutes with filtered inductive limits.
\end{theorem}

\bpr
From the theory of categories, $F$ has a unique extension to a functor 
$\tilde F\colon \Fun(\FMod^\op,\Sets)\to \Mod$ 
commuting with all inductive limits (\cite{Sheaves}, Corollary I.5.4). 
Here $\FMod$ embeds as usual into $\Fun(\FMod^\op,\Sets)$ through the Yoneda embedding 
$M\mapsto \Hom(-,M)$, which, in turn, has a canonical extension to a functor 
$J\colon \Mod\to\Fun(\FMod^\op,\Sets)$, defined by the self-same formula. 
Since limits of functors are calculated point-wise 
and $\Hom(P,-)$, for $P$ free and finitely generated, is right-exact and commutes with filtered 
inductive limits, the same is true of the composite $\hat F = \tilde F \circ J$. 
Uniqueness of $\hat F$ is evident. 
\epr

\section{Polynomial and Strict Polynomial Functors}

Let us now bring to mind the classical notions of polynomial maps and functors 
(Eilenberg \& MacLane \cite{EM}, sections 8--9).  
Let $\phi\colon M\to N$ be a (non-linear) map of modules. The $n$'th 
\textbf{deviation} of $\phi$ is the map 
$$ 
\dev{\phi}{x_1}{x_{n+1}} = \sum_{I\subseteq[n+1]} (-1)^{n+1-\abs{I}} \phi\left(\sum_{i\in I} x_i\right) 
$$
in $n+1$ variables. 
The map $\phi\colon M\to N$ is \textbf{polynomial} of degree $n$ if its $n$'th deviation vanishes,  
and the functor $F\colon \FMod \to \Mod$ is \textbf{polynomial} of degree $n$ 
if every arrow map $\Hom(M,N) \to \Hom(F(M),F(N))$ is. 

We next recall the strict polynomial maps (\emph{lois polynomes}) 
from the work of Roby (\cite{Roby}, paragraph 1.2)
and the strict polynomial functors introduced by Friedlander \& Suslin (\cite{FS}, Definition 2.1). 
A \textbf{strict polynomial map} is a natural transformation 
$\phi\colon M\otimes - \to N\otimes -$ between functors ${}_\B\CAlg\to \Sets$, 
where ${}_\B\CAlg$ designates the category of commutative, unital $\B$-algebras, 
and $\Sets$ is the category of sets. 
The functor $F\colon \FMod \to \Mod$ is  \textbf{strict polynomial} of degree $n$ if the arrow maps
$\Hom(M,N) \to \Hom(F(M),F(N)) $
have been given a (multiplicative) strict polynomial structure. 

Strict polynomial maps and functors decompose as the direct sum of their homogeneous components 
(\cite{Roby}, Proposition I.4). 
We denote by $\HPol_n$ the abelian category of homogeneous functors of degree $n$, where a
\textbf{natural transformation} $\eta\colon F\to G$ of homogeneous functors is a family of homomorphisms
$\eta_M\colon F(M)\to G(M)$, for $M\in\FMod$,  
such that for any modules $M$ and $N$, any algebra $A$ and any $\omega\in A\otimes \Hom(M,N)$, 
there is a commutative diagram: 
$$ 
\xymatrix{
A\otimes F(M) \ar[d]_{F(\omega)} \ar[r]^{1\otimes \eta_M} & A\otimes G(M) \ar[d]^{G(\omega)} \\
A\otimes F(N) \ar[r]_{1\otimes \eta_N} & A\otimes G(N) 
} 
$$

Strict polynomial functors are not determined by their underlying functors 
since the strict structure supplies 
auxiliary data. The example below should serve as a warning. (In fact, this failure occurs 
already at the level of strict polynomial maps; cf.~Example 7 of \cite{PM}.)

\bex 					\label{X: Frobenius} 
Let $\B=\Z$. The functor $F=\Z/p\Z\otimes-$, for $p$ a prime, is an homogeneous linear functor. 
Yet, it can also be equipped with  
an homogeneous structure of degree $p$, viz. the \emph{Frobenius twist}
\begin{gather*}
F^{(1)}\colon A\otimes \Hom(M,N) \to A\otimes \Hom(\Z/p\Z\otimes M,\Z/p\Z\otimes N) \\
\sum_i a_i \otimes \alpha_i \mapsto \sum_i a_i^p\otimes \alpha_i. \qedhere
\end{gather*}
\eex

Let $\Gamma^n$ denote the functor of divided $n$'th powers. For the matrix ring~$\B^{n\times n}$, 
the divided power module $\GB$ becomes an algebra under the 
multiplication $\alpha^{[n]}\circ \beta^{[n]} = (\alpha\beta)^{[n]}$. 
(This operation is not to be confused with the ever-present graded multiplication on any 
divided power algebra $\Gamma(M)$.)

That strict polynomial functors arise as modules over the Schur algebra 
was first proved by Friedlander \& Suslin for finite fields (\cite{FS}, Theorem 3.1), 
and later in full generality in co-operation with Bendel (\cite{SFB}, Theorem 2.4): 

\begin{theorem}			\label{S: HPol} 
The functor $\GH$ is a 
small projective generator for $\HPol_n$, through which there is an equivalence of categories
$$
\Psi\colon\HPol_n \to \GBMod, \qquad F \mapsto F(\B^n).
$$ 
The left module structure on $F(\B^n)$ is given by $\tau^{[n]}\cdot x = F(\tau)(x)$. 
\end{theorem}

\section{Numerical Functors}				\label{A: Definitions}

From now on, we shall assume our base ring $\B$ to be \emph{binomial}
in the sense of Hall \cite{Hall}; that is, commutative, unital and torsion-free, so that it naturally 
embeds into $\Q\otimes_\Z \B$, and closed under the formation of \emph{binomial co-efficients} 
$$
r\mapsto \binom{r}{n} = \frac{r(r-1)\cdots(r-n+1)}{n!}, \qquad n\in \B.
$$    
Examples include $\Z$, the ring $\hat{\Z}_p$ of $p$-adic integers, as well as all $\Q$-algebras. 
These rings were extensively studied in \cite{Elliott} and \cite{BR}.

Following \cite{PM} (Definition 5), we say that 
the polynomial map $\phi\colon M\to N$ is \textbf{numerical} of degree $n$ if it satisfies the equation 
$$
\phi(rx) = \sum_{k=0}^n \binom{r}{k}\Dev{\phi}{_k x}, \qquad r\in\B,\ x\in M.
$$

We propose the following definition, extending Eilenberg and MacLane's notion of 
polynomial functors. 

\bdf 			\label{D: Numerical} 
The functor $F\colon \FMod \to \Mod$ is \textbf{numerical} of degree $n$ if every arrow map 
$F\colon \Hom(M,N) \to \Hom(F(M),F(N))$ is. 
The category of numerical functors of degree $n$ will be denoted 
by the symbol $\Num_n$. 
\edf 

So, for example, a functor is numerical of degree $0$ if and only if it is constant. It is 
of degree $1$ if and only if it is the translate of a linear functor. 

By simple algebraical considerations, $\Num_n$ 
may be verified to be abelian (the case $\B=\Z$ is well known). 
It is, moreover, closed under direct sums. We shall presently see 
that it possesses a small projective generator.

Every strict polynomial functor is numerical, and if the base ring $\B$ is 
a $\Q$-algebra, the two strains co-incide. 
Over $\Z$, the notions of polynomial and numerical functor may be equated, for then 
all polynomial maps are automatically numerical. 
These assertions are consequences of the corresponding statements for maps. 
Confer the remarks succeeding Definition 5 in \cite{PM}.


An alternative description of numerical functors may be given, matching 
 the definition given of strict polynomial functors.  
Suppose that, for each finitely generated and free module $M$, a module $F(M)$ is given. 
Moreover, for any two finitely generated, free modules $M$ and $N$, suppose given maps  
$$
F_A \colon A\otimes \Hom(M,N) \to A\otimes \Hom(F(M),F(N)),
$$
natural in the binomial algebra $A$. We say that these maps $F_A$ are \emph{multiplicative} if
$F_A(1_A\otimes 1_M) = 1_A\otimes 1_{F(M)}$ and
$$
F_A\left(\sum_{i=1}^p a_i\otimes \alpha_i\right) F_A\left(\sum_{j=1}^q b_j\otimes \beta_j\right) 
= F_A\left(\sum_{i=1}^p\sum_{j=1}^q a_ib_j \otimes \alpha_i\beta_j\right)
$$ 
for homomorphisms $\alpha_i\colon N\to P$ and $\beta_j\colon M\to N$.

\bth				\label{S: BAlg-polynomial}
The functor $F\colon\FMod\to\Mod$ is numerical of degree $n$ if and only if its arrow maps extend to a system  
$$
F_A \colon A\otimes \Hom(M,N) \to A\otimes \Hom(F(M),F(N))
$$
of maps of degree $n$, multiplicative and natural in the binomial algebra $A$.
\eth

\bpr
According to Theorem 10 of \cite{PM}, the module map $F\colon \Hom(M,N)\to\Hom(F(M),F(N))$ 
is numerical of degree $n$ precisely when it extends to a degree $n$ natural transformation 
$$
F_A \colon A\otimes \Hom(M,N) \to A\otimes \Hom(F(M),F(N)),
$$
given by the formula 
$$
F_A(a_1\otimes\alpha_1+\cdots a_k\otimes \alpha_k) 
= \sum_{m_1,\dots,m_k=0}^\infty \binom{a_1}{m_1}\cdots\binom{a_k}{m_k} 
\otimes F\left(\De_{m_1}\alpha_1 \de \cdots \de \De_{m_k}\alpha_k\right).
$$
Only a finite number of summands are non-zero by the polynomiality of $F$.
Since this assignment is functorial, multiplicativity of the system is clear.
\epr

Natural transformations may also be given a more involved rendition, 
comparable to transformations of strict polynomial functors:

\bth 
Let $\eta\colon F\to G$ be a natural transformation of numerical functors. 
For any modules $M$ and $N$, any binomial algebra $A$ and any $\omega\in A\otimes \Hom(M,N)$, 
the following diagram commutes: 
\beq					\label{E: Natural transformation}
\xymatrix{
A\otimes F(M) \ar[d]_{F(\omega)} \ar[r]^{1\otimes \eta_M} & A\otimes G(M) \ar[d]^{G(\omega)} \\
A\otimes F(N) \ar[r]_{1\otimes \eta_N} & A\otimes G(N) 
} 
\eeq
\eth

\bpr 
Consider homomorphisms $\alpha_1,\dots,\alpha_k\colon M\to N$. By Theorem 9 of \cite{PM}, there exist 
homomorphisms $\beta_{m_1,\dots,m_k}$ and $\gamma_{m_1,\dots,m_k}$ 
(the indices ranging over natural numbers, 
and only a finite number being non-zero), such that 
\begin{align*}
F(a_1\otimes \alpha_1 + \cdots + a_k\otimes\alpha_k) 
&= \sum_{m_1,\dots,m_k} \binom{a_1}{m_1}\cdots \binom{a_k}{m_k} \otimes \beta_{m_1,\dots,m_k} \\
G(a_1\otimes \alpha_1 + \cdots + a_k\otimes\alpha_k) 
&= \sum_{m_1,\dots,m_k} \binom{a_1}{m_1}\cdots \binom{a_k}{m_k} \otimes \gamma_{m_1,\dots,m_k} 
\end{align*}
for any $a_1,\dots,a_k$ in any binomial algebra $A$.
The naturality of $\eta$ ensures that 
$$ 
\sum_{m_1,\dots,m_k} \binom{a_1}{m_1}\cdots \binom{a_k}{m_k} \eta_N\beta_{m_1,\dots,m_k} 
= \sum_{m_1,\dots,m_k} \binom{a_1}{m_1}\cdots \binom{a_k}{m_k}  \gamma_{m_1,\dots,m_k}\eta_M .
$$
Specialise first to the case $a_2=a_3=\dots=0$, to obtain 
$$ 
\sum_{m_1} \binom{a_1}{m_1} \eta_N\beta_{m_1,0,\dots} 
= \sum_{m_1} \binom{a_1}{m_1} \gamma_{m_1,0,\dots}\eta_M .
$$
Successively putting $a_1=0,1,2,\dots$ leads to 
$\eta_N\beta_{m_1,0,\dots} = \gamma_{m_1,0,\dots}\eta_M$
for all $m_1$. Proceeding inductively, one shews that 
$\eta_N \beta_{m_1,\dots,m_k} = \gamma_{m_1,\dots,m_k}\eta_M$
for any assortment of indices. 

Commutativity of the diagram \eref{E: Natural transformation} 
for $\omega=a_1\otimes \alpha_1 + \cdots + a_k\otimes\alpha_k$ is then easily demonstrated:
\begin{multline*} 
(1\otimes\eta_N)F(\omega) 
= (1\otimes\eta_N)\left( \sum_{m_1,\dots,m_k} \binom{a_1}{m_1}
\cdots \binom{a_k}{m_k} \otimes \beta_{m_1,\dots,m_k}\right) \\
= \sum_{m_1,\dots,m_k} \binom{a_1}{m_1}\cdots \binom{a_k}{m_k} \otimes \eta_N\beta_{m_1,\dots,m_k}
= \sum_{m_1,\dots,m_k} \binom{a_1}{m_1}\cdots \binom{a_k}{m_k} \otimes \gamma_{m_1,\dots,m_k}\eta_M \\
= \left( \sum_{m_1,\dots,m_k} \binom{a_1}{m_1}\cdots \binom{a_k}{m_k} 
\otimes \gamma_{m_1,\dots,m_k}\right) (1\otimes\eta_M) 
= G(\omega)(1\otimes\eta_M). \qedhere
\end{multline*}
\epr

\section{Properties of Numerical Functors}

We give some equivalent characterisations of numericality, which may perhaps be more convenient in practice. 

\bth 			
The following conditions, assumed to hold for all scalars $r$ and all homomorphisms $\alpha$, 
are equivalent to place on a polynomial functor $F$ of degree $n$. 
\balph 
\item[\upshape{A.}] 
$\displaystyle F(r\alpha) = \sum_{k=0}^n \binom{r}{k}\Dev{F}{_k \alpha}$ 
(the definition of numerical functor)
\item[\upshape{A$'$.}] 
$\displaystyle F(r 1_{\B^n}) = \sum_{k=0}^n \binom{r}{k}\Dev{F}{_k 1_{\B^n}}$
\item[\upshape{B.}] 
$\displaystyle F(r\alpha) = \sum_{m=0}^n (-1)^{n-m}\binom{r}{m}\binom{r-m-1}{n-m} F(m\alpha)$
\item[\upshape{B$'$.}] 
$\displaystyle F(r 1_{\B^n}) = \sum_{m=0}^n (-1)^{n-m}\binom{r}{m}\binom{r-m-1}{n-m} F(m 1_{\B^n})$ 
\ealph
\eth

\bpr 
That A and B are equivalent follows from Theorem 7 of \cite{PM}, 
as does the equivalence of A$'$ and B$'$. Clearly B implies B$'$, so there 
remains to establish the implication of B by B$'$. Hence assume B$'$, and put 
$$ 
Z_m = (-1)^{n-m}\binom{r}{m}\binom{r-m-1}{n-m}. 
$$
In so far as $1_{\B^q}$ factorises through $1_{\B^n}$, which holds in case $q\leq n$, the equation 
\beq 			\label{E: Auxiliary numerical}
F(r 1_{\B^q}) = \sum_{m=0}^n Z_m F(m 1_{\B^q})
\eeq
will clearly hold by B$'$. Now assume \eref{E: Auxiliary numerical} 
holds for a certain $q\geq n$. 
Let $\pi_i\colon\B^{q+1}\to\B^{q+1}$ ($1\leq i\leq q+1$) denote the canonical projections. 
Using the fact that 
$$
0 = F(\alpha_1\de\cdots\de\alpha_{q+1}) 
= \sum_{I\subseteq[q+1]} (-1)^{q+1-\abs{I}} F\left(\sum_{i\in I} \alpha_i\right)
$$
(the polynomiality of $F$), we compute:
\begin{multline*}
F(r1_{\B^{q+1}}) = F(r\pi_1+\cdots+r\pi_{q+1}) 
= \sum_{I\subset[q+1]} (-1)^{q-\abs{I}} F\left(\sum_{i\in I} r\pi_i\right) \\
= \sum_{I\subset[q+1]} (-1)^{q-\abs{I}} \sum_{m=0}^n Z_m F\left(\sum_{i\in I} m\pi_i\right) 
= \sum_{m=0}^n Z_m \sum_{I\subset[q+1]} (-1)^{q-\abs{I}}  F\left(\sum_{i\in I} m\pi_i\right) \\
= \sum_{m=0}^n Z_m F(m\pi_1+\cdots +m\pi_{q+1} ) = \sum_{m=0}^n Z_m F(m 1_{\B^{q+1}} ),
\end{multline*}
so that \eref{E: Auxiliary numerical} holds for all $q$ by induction. 
Finally, in the case of an arbitrary homomorphism $\alpha\colon \B^p\to \B^q$, we have
$$
F(r\alpha) = F(r1_{\B^q})F(\alpha) 
 = \sum_{m=0}^n Z_m F(m1_{\B^q} )F(\alpha) = \sum_{m=0}^n Z_m F(m\alpha),
$$
and we have proved B. 
\epr

The following very pleasant formula is an immediate consequence of the corresponding 
formula for maps, Theorem 8 in \cite{PM}. 

\bth 				\label{S: Pleasant}
The module functor $F$ is numerical of degree $n$ if and only if, for any scalars $a_i$ and homomorphisms $\alpha_i$, the following equation holds: 
$$ 
\dev{F}{a_1\alpha_1}{a_k\alpha_k} 
= \sum_{\substack{g_1+\cdots+g_k\leq n \\ g_i\geq 1}} \binom{a_1}{g_1}\cdots\binom{a_k}{g_k} 
F\left(\De_{g_1}\alpha_1 \de \cdots \de \De_{g_k} \alpha_k\right). 
$$ 
\eth


\section{Quasi-Homogeneous Functors}			\label{A: Quasi-homogeneous}

A module map $\phi\colon M\to N$ 
naturally extends to a map $\phi\colon \B[M]\to N$, where $\B[M]$ denotes the free module on $M$. The 
numericality of $\phi$ may then be reformulated by demanding that the extended map 
$\phi\colon \B[M]\to N$ factorise through the canonical homomorphism $\delta_n\colon M\to \BM=\B[M]/I$, 
where 
$$		
I =  \Ideal{ [x_1\de\cdots\de x_{n+1}] |  x_i\in M }  
 + \Ideal{ [rx] - \sum_{k=0}^n \binom{r}{k}\left[\De_k x\right] | 
 r\in \B, \ x\in M }. 
$$
We say that the numerical map $\phi$ is \textbf{quasi-homogeneous} if it vanishes on 
 $\BM\cap J$, 
where $J$ is the $\Q$-linear subspace generated in $\Q\otimes_\Z \BM$ by all elements of the form 
$[rx] - r^n[x]$, for $r\in \B$ and $x\in M$. 
 
\bdf 			\label{D: Quasi-homogeneous} 
The numerical functor $F$ is \textbf{quasi-homogeneous} 
if every arrow map $\Hom(M,N) \to \Hom(F(M),F(N))$ is.
\edf

Being quasi-homogeneous is a necessary condition for a functor to 
admit a (strict polynomial) homogeneous structure. We shall later give 
a sufficient condition. 

One may enquire why the simpler condition $F(r\alpha)=r^n F(\alpha)$ does not suffice.
While this would serve the purpose equally well in degrees $2$ and $3$, it will be inadequate in 
degree $4$. See Example \ref{X: Quartical} below. 

\bth
The functor $F$, numerical of degree $n$, 
is quasi-homogeneous of degree $n$ if and only if it annihilates $P_n(\Hom(M,N))\cap K$ 
for all modules $M$, $N$, where $K$ is the $\Q$-linear subspace generated in 
$\Q\otimes_\Z P_n(\Hom(M,N))$ by all elements of the form 
$$
[\alpha_1\de\cdots\de\alpha_k] - \sum_{\substack{g_1+\cdots+g_k=n \\ g_i\geq 1}} 
\frac{1}{g_1!\cdots g_k!}
\left[ \De_{g_1}\alpha_1 \de \cdots \de \De_{g_k} \alpha_k \right].
$$
\eth

\bpr
The special case $k=0$ of the above relation gives $[0]=[\de]\equiv 0 \mod K$.
Now, using Theorem \ref{S: Pleasant},
$$
[r\alpha] \equiv [\de r\alpha] \equiv \frac{1}{n!} \left[ \De_n r\alpha \right] 
= \frac{r^n}{n!} \left[ \De_n \alpha \right] \equiv r^n[\de\alpha] \equiv r^n [\alpha] \mod K.
$$
On the other hand, again using Theorem \ref{S: Pleasant}, 
\begin{multline*}
r^n [\alpha_1\de\cdots\de\alpha_k] \equiv [r\alpha_1\de\cdots\de r\alpha_k] \\
= \sum_{\substack{g_1+\cdots+g_k=n \\ g_i\geq 1}} 
\binom{r}{g_1} \cdots \binom{r}{g_k}
\left[ \De_{g_1}\alpha_1 \de \cdots \de \De_{g_k} \alpha_k \right] \mod J.
\end{multline*}
Since $\Q\otimes_\Z P_n(\Hom(M,N))$ is torsion-free, one may identify the co-efficients of $r^n$, 
whence
$$
[\alpha_1\de\cdots\de\alpha_k] \equiv 
\sum_{\substack{g_1+\cdots+g_k=n \\ g_i\geq 1}} 
\frac{1}{g_1!\cdots g_k!}
\left[ \De_{g_1}\alpha_1 \de \cdots \de \De_{g_k} \alpha_k \right] \mod J. \qedhere
$$
\epr

\bex
A quadratic functor $F$ is easily verified to be quasi\hyp homogeneous if and only if 
$$
F(\de) = 0 \qquad\text{and}\qquad 2 F(\de \alpha) = F(\alpha\de\alpha),  
$$
or if and only if $F(r\alpha)=r^2 F(\alpha)$. 

In the cubical case, the conditions are 
\begin{gather*}
F(\de) = 0, \qquad 6 F(\de\alpha) = F(\alpha\de\alpha\de\alpha) 
\qquad\text{and} \\
 2 F(\alpha\de\beta) = F(\alpha\de\alpha\de\beta) + F(\alpha\de\beta\de\beta),
\end{gather*}
or simply $F(r\alpha)=r^3 F(\alpha)$.
\eex

\bex				\label{X: Quartical}
Let $\B=\Z$, and consider the functor 
$F(M)=\Z/12 \otimes \Gamma^2(M)$, which is quadratic, \emph{a fortiori} quartical. 
The reader might be led into imagining $F$ to be quasi-homogeneous of degree $4$, since it satisfies 
$$
F(r\alpha) = 1 \otimes (r\alpha)^{[2]} = r^2 \otimes \alpha^{[2]} 
= r^4 \otimes \alpha^{[2]} = r^4 F(\alpha).
$$
The abelian group generated in $P_4(\Hom(M,N))$ by all elements of the form 
$[r\alpha] - r^4[\alpha]$ contains, 
for any homomorphisms $\alpha$ and $\beta$, the element 
$$
12[\alpha\de\beta] 
- 2[\alpha\de\alpha\de\beta] - 2[\alpha\de\alpha\de\beta\de\beta] - 2[\alpha\de\beta\de\beta\de\beta],
$$
which $F$ will indeed annihilate. Dividing by $2$, however, produces the element
$$
6[\alpha\de\beta] 
- [\alpha\de\alpha\de\beta] - [\alpha\de\alpha\de\beta\de\beta] - [\alpha\de\beta\de\beta\de\beta]
$$
in $J$, which is not annihilated, so that $F$ is, in fact, not quasi-homogeneous.
\eex

\section{The Morita Equivalence}

We now exhibit a projective generator of the category of numerical functors, 
demonstrating the equivalence of $\Num_n$ with a suitable module category. 

Throughout, let $K$ be a fixed free and finitely generated module. 
The composite functor $\BHK$, transforming $\chi\colon M\to N$ unto
$$
[\chi_\ast]\colon \BHK[M] \to \BHK[N],
$$
is plainly numerical of degree $n$. 

%

\bth[The Yoneda Lemma] 
Let $F$ be numerical of degree $n$. The map
$$
\Upsilon_{K,F}\colon \Nat( \BHK[-], F) \to F(K), \qquad \eta \mapsto \eta_K([1_K]),
$$
is an isomorphism of modules. It is natural 
in the sense that the following two diagrams commute: 
\begin{gather*} 
\xymatrix{
K \ar[d]_{\beta} & \Nat(\BHK,F) \ar[d]_{[\beta^\ast]^\ast} \ar[r]^-{\Upsilon_{K,F}} 
& F(K) \ar[d]^{F(\beta)} \\
L & \Nat(P_n\Hom(L,-),F) \ar[r]_-{\Upsilon_{L,F}} & F(L)  
}
\\
\xymatrix{
F \ar[d]_{\xi} & \Nat(\BHK,F) \ar[d]_{\xi_\ast} \ar[r]^-{\Upsilon_{K,F}} & F(K) \ar[d]^{\xi_K } \\
G & \Nat(\BHK,G) \ar[r]_-{\Upsilon_{K,G}} & G(K)  
}
\end{gather*} 
\eth

\bpr
The proof is the usual one. Naturality is obvious. 
Consider the following commutative diagram:
$$ 
\xymatrixcolsep{1.2pc} 
\xymatrix{
K \ar[d]_\alpha & \BHK[K] \ar[d]_{[\alpha_\ast]} \ar[r]^-{\eta_K} & F(K) \ar[d]^{F(\alpha)} 
& [1_K] \ar[d] \ar[r] & \eta_K([1_K]) \ar[d] \\
M & \BHK[M] \ar[r]_-{\eta_M} & F(M) & [\alpha] \ar[r] & \eta_M([\alpha]) = F(\alpha)(\eta_K([1_K]))
} 
$$
Upon inspection, we find that $\Upsilon_{K,F}$ has an inverse transforming $y\in F(K)$ to
$$ 
\eta_M\colon \BHK[M] \to F(M) , \qquad [\alpha] \mapsto F(\alpha)(y).
$$
The numericality of $F$ is used in a most essential way to ensure that the map 
$\Hom(K,M)\to \Hom(F(K),F(M))$ factorise through $\BHK[M]$.
\epr

Putting $F=\BHK$, we obtain a module isomorphism
$$
\Upsilon\colon \Nat( \BHK[-], \BHK[-] ) \to P_n(\End K), 
$$
whose inverse transforms $[\sigma]$, for $\sigma\colon K\to K$, to the natural transformation 
$$
[\sigma^\ast]\colon  \BHK[-] \to \BHK[-].
$$ 

The module $\B[\End K]$ comes equipped with a 
 multiplication, given by $[\sigma]\circ [\tau]=[\sigma\tau]$, 
which descends unto the algebras $P_n(\End K)$. 
(This is not to be confused with another multiplication, defined on any $\B[M]$, 
given by $[x][y]=[x+y]$.) Under the Yoneda correspondence, 
$$
\Upsilon^{-1}([\sigma]\circ [\tau]) = \Upsilon^{-1}([\sigma\tau]) =  [(\sigma\tau)^\ast] \\
= [\tau^\ast] \circ [\sigma^\ast] = \Upsilon^{-1}([\tau])\circ \Upsilon^{-1}([\sigma]),
$$
so that the multiplications are reversed by $\Upsilon$:

\bth 			\label{S: Anti-iso}
The Yoneda correspondence provides an anti-isomorphism of rings
$$ 
\Nat(\BHK[-],\BHK[-]) \cong P_n(\End K). 
$$
\eth

We now prove that $\BH$ provides a projective generator for $\Num_n$.

\blem 			\label{L: Vanishing Lemma Num} 
A polynomial functor of degree $n$ that vanishes on $\B^n$ is identically zero. 
\elem

\bpr 
Let $F$ be such a functor. We shall shew that $F(\B^q)=0$ for all natural numbers $q$.
When $q\leq n$, then $\B^q$ is a direct summand of $\B^n$, and therefore 
$F(\B^q)$ will be a direct summand of $F(\B^n)=0$. 

Proceeding by induction, suppose $F(\B^q)=0$ for some $q\geq n$. We shew that also $F(\B^{q+1})=0$.
Let $\pi_i\colon \B^q\to \B^q$ ($1\leq i\leq q+1$) denote the canonical projections. 
Since $F$ is polynomial of degree $n$, we have 
$$
0= F(\pi_1\de\cdots\de\pi_{q+1}) 
=\sum_{I\subseteq [q+1]} (-1)^{q+1-\abs{I}} F\left(\sum_{i\in I} \pi_i\right). 
$$
Consider an $I$ with $\abs{I}\leq q$. Since $\sum_{I} \pi_i$ factorises through 
$\B^q$, the homomorphism $F\left(\sum_{I} \pi_i\right)$ will factorise through 
$F(\B^q)=0$.  Only $I=[q+1]$ will give a non-trivial contribution to the sum above, yielding 
$$
0=F(\pi_1 + \cdots + \pi_{q+1}) = F(1_{\B^{q+1}}) = 1_{F(\B^{q+1})}. \qedhere
$$
\epr

\bth 				\label{S: Num} 
The functor $\BH$ is a small projective generator for $\Num_n$, through 
which there is an equivalence of categories
$$ 
\Phi\colon \Num_n \to \BBMod, \qquad F\mapsto F(\B^n).
$$ 
The left module structure on $F(\B^n)$ is given by $ [\tau]\cdot x = F(\tau)(x)$. 
\eth

\bpr 
Write $Q_n = \BH$. We first establish the projectivity of $Q_n$ by shewing that $\Nat(Q_n,-)$ 
is right-exact, i.e.~preserves epimorphisms. If $\eta\colon F\to G$ is epimorphic, then 
 $\eta_{\B^n}\colon F(\B^n) \to G(\B^n)$ is surjective, and corresponds to 
$\eta_\ast\colon \Nat(Q_n,F) \to \Nat(Q_n,G)$ 
under the Yoneda map.

Since $0=\Nat(Q_n,F)\cong F(\B^n)$ implies $F=0$ by the lemma, the functor $Q_n$ is a generator. 

For a direct sum $\bigoplus F_k$ of functors, we compute, using the Yoneda Lemma,
$$
\Nat\left(Q_n,\bigoplus F_k\right) \cong \left(\bigoplus F_k\right)(\B^n) = \bigoplus F_k(\B^n) 
\cong \bigoplus \Nat\big(Q_n,F_k\big) , 
$$
establishing that $\Nat(Q_n,-)$ preserves direct sums, i.e.~smallness of $Q_n$. 

As $\Num_n$ is an abelian category with arbitrary direct sums, a Morita equivalence will arise through 
the projective generator:
$$ 
\xymatrixcolsep{3pc}
\xymatrix{ 
\Num_n \upar[r]^-{\Nat(Q_n,-)} & \upar[l]^-{Q_n \otimes_S -} \SMod 
} 
$$
The underlying ring for the module category is, by Theorem \ref{S: Anti-iso},
$$
S=\Nat(Q_n,Q_n)^\op\cong P_n(\End \B^n) =\BB.
$$ 

A functor $F$ corresponds to the abelian group $\Nat(Q_n,F) \cong F(\B^n)$.
Under the Yoneda map, an element $x\in F(\B^n)$ will correspond to the natural transformation
$$ 
\eta_M\colon Q_n(M) \to F(M), \qquad [\alpha] \mapsto F(\alpha)(x). 
$$
In like wise, $[\tau]\in \BB$ will correspond to 
$$ 
\theta_M\colon Q_n(M) \to Q_n(M), \qquad [\alpha] \mapsto [\alpha\circ \tau]. 
$$
The product of $\theta$ and $\eta$ is the transformation 
$$ 
(\eta\circ\theta)_M\colon  Q_n(M) \to F(M), \qquad [\alpha] \mapsto F(\alpha\circ \tau)(x), 
$$
which under the Yoneda map corresponds to 
$$
(\eta\circ\theta)_{\B^n}([1_{\B^n}]) = F(1_{\B^n} \circ \tau)(x) = F(\tau)(x) \in F(\B^n),
$$
affirming the module structure on $F(\B^n)$. 
\epr

\section{The Divided Power Map} 

Taking as our starting point a finitely generated and free 
module $M$, we next propose a study of the \textbf{divided power map}
$$
\gamma_n \colon \BM \to \GM, \qquad [x] \mapsto x^{[n]}.
$$
It is linear, being induced by the strict polynomial, \emph{a fortiori} numerical, map
$M \to \GM$, given by $x\mapsto x^{[n]}$.

\blem 			\label{L: gamma_n}
If $x_1,\dots,x_k\in M$, then 
$$
\gamma_n([x_1 \de \cdots \de x_k]) 
= \sum_{\substack{g_1+\cdots+g_k=n \\ g_i\geq 1}} x_1^{[g_1]}\cdots x_k^{[g_k]} . 
$$ 
\elem 

\bpr 
By the definition of deviations, 
$$
\gamma_n([x_1 \de \cdots \de x_k]) 
= \sum_{I\subseteq [k]} (-1)^{k-\abs{I}} \left(\sum_{i\in I} x_i\right)^{[n]} .
$$  
The co-efficient of a given monomial $x_1^{[g_1]}\cdots x_k^{[g_k]}$ ($g_1+\cdots+g_k=n$) will be
$$ 
\sum_{\set{i \, | \, g_i>0} \subseteq I\subseteq [k]} (-1)^{k-\abs{I}} = 
\begin{cases}
1 & \text{if $\set{i \, | \, g_i>0}=[k]$,} \\
0 & \text{otherwise.}
\end{cases}  
\qedhere
$$
\epr

%
%

In general, the divided power map is neither injective nor surjective:

\bth 				\label{S: Kernel}
The divided power map $\gamma_n\colon \BM\to \GM$ has kernel
$$ 
\Ker \gamma_n = \Q\otimes_{\Z} \gen{[rz] - r^n[z] | r\in\B, \ z\in M } \cap \BM. 
$$ 
If the rank of $M$ is at least $n$, it 
is surjective precisely when $(n-1)!$ is inversible in $\B$.
\eth

\bpr 
We begin with investigating the kernel. Denoting 
$$ 
L = \Q\otimes_{\Z} \gen{[rz] - r^n[z] | r\in\B, \ z\in M } ,
$$
we evidently have $L\cap\BM\subseteq \Ker\gamma_n$. 

For the reverse inclusion, let $z\in M$ and $m\in\N$. Calculating modulo $L$:
\begin{multline} 			\label{E: Modulo L}
\left[\De_m z\right] 
= \sum_{I\subseteq [m]} (-1)^{m-\abs{I}} \left[ \sum_{i\in I} z\right] 
= \sum_{j=0}^m (-1)^{m-j} \binom mj [ jz]  \\
\equiv \sum_{j=0}^m (-1)^{m-j} \binom mj j^n [ z] 
= m!\stir{n}{m} [z], 
\end{multline}
where $\stir{n}{m}$ denotes a Stirling number of the second kind.

Let $\{e_1,\dots,e_k\}$ be a basis for $M$. A basis for $\BM$ is then, by 
Theorem 4 of \cite{PM}, constituted by the elements
$$ 
[e_{p_1}\de\cdots\de e_{p_m}], 
\qquad 1\leq p_1\leq\dots\leq p_m\leq k, \quad 0\leq m\leq n.
$$
Letting $[f_1\de\cdots\de f_m]$ be one such basis element of $P_n(M)$, we  have, 
by \eref{E: Modulo L}:
\begin{multline*}
m!\stir{n}{m} [f_1\de\cdots\de f_m] 
=  m!\stir{n}{m} \sum_{I\subseteq [m]} (-1)^{m-\abs{I}} \left[ \sum_{i\in I} f_i\right] 
\equiv \sum_{I\subseteq [m]} (-1)^{m-\abs{I}} \left[ \De_m \sum_{i\in I} f_i\right] \\
= \sum_{I\subseteq [m]} (-1)^{m-\abs{I}} \sum_{\emptyset\subset I_1,\dots,I_m\subseteq I} 
\left[ \De_{i_1\in I_1} f_{i_1} \de \cdots \de \De_{i_m\in I_m} f_{i_m} \right] \mod L. 
\end{multline*}
We may write this as $\xi+\xi'$, where $\xi$ is a sum of deviations of order exactly~$m$, 
and $\xi'$ collects the higher-order deviations. 
An $m$'th-order deviation arises precisely when the sets $I_1,\dots,I_m$ are singletons. 
Supposing $g_i$ of the sets $I_1,\dots,I_m$ are equal to $\{i\}$, for $i=1,\dots,m$, 
we may express $\xi$ as
\begin{multline*}
\xi = 
\sum_{I\subseteq [m]} (-1)^{m-\abs{I}} \sum_{\substack{g_1+\cdots+g_m=m \\ g_i>0 \ek i\in I }} 
\binom{m}{g_1,\dots,g_m} \left[ \De_{g_1} f_1 \de \cdots \de \De_{g_m} f_m \right] \\
= 
\sum_{\substack{g_1+\cdots+g_m=m \\ g_i\geq 0 }} 
\left( \sum_{\set{i \, | \, g_i>0}\subseteq I\subseteq [m]} (-1)^{m-\abs{I}} \right)  
\binom{m}{g_1,\dots,g_m} \left[ \De_{g_1} f_1 \de \cdots \de \De_{g_m} f_m \right].
\end{multline*}
The inner sum is $1$ when all $g_i=1$, and $0$ otherwise, leading to
$$
\xi = 
\binom m{1,\dots,1} [f_1\de\cdots\de f_m] = m! [f_1\de\cdots\de f_m].
$$
We thus have 
$$ 
m! \stir{n}{m} [f_1\de\cdots\de f_m]  \equiv m! [f_1\de\cdots\de f_m] + \xi'  \mod L,
$$
and, consequently, provided that $1<m<n$ (so that $\stir{n}{m} > 1$),
$$ 
[f_1\de\cdots\de f_m] \equiv  \frac {1}{m! \left( \stir{n}{m} - 1\right)} \cdot \xi'  \mod L. 
$$

In summary, what we have proved is that 
a deviation of order $m$ may be expressed, modulo $L$, in terms of deviations of 
higher orders, provided $1<m<n$. It may happen that $m=0$ or $m=1$. We may then use the relations 
$$ 
[\de] = [0] \equiv 0 
\qquad\text{and}\qquad
[\de f] = [f]-[0] \equiv \frac{1}{n! \stir{n}{n} } \left[\De_n f\right] 
= \frac{1}{n!} \left[\De_n f\right] \mod L, 
$$
which are both consequences of \eref{E: Modulo L}, to achieve the same end.

Now suppose $\omega\in\Ker\gamma_n$. Successively using the above relations, 
we may express $\omega$ as a fractional linear combination 
$$ 
\omega \equiv 
\sum_{\substack{g_1+\cdots+g_k=n \\ g_i\geq 0}} c_{g_1,\dots,g_k}
\left[\De_{g_1} e_1 \de \cdots \de \De_{g_k} e_k \right]  \mod L, 
\qquad c_{g_1,\dots,g_k}\in \Q\otimes_\Z \B,
$$ 
of deviations of rank $n$ of the basis elements $e_i$. Apply $\gamma_n$ and use the lemma: 
$$ 
0=\gamma_n(\omega) 
= \sum_{\substack{g_1+\cdots+g_k=n \\ g_i\geq 0}} c_{g_1,\dots,g_k}
e_1^{g_1}\cdots e_k^{g_k}. 
$$ 
Because the elements $e_1^{g_1}\cdots e_k^{g_k}$ 
are independent in $\Gamma^n(M)$, it must be that all co-efficients $c_{g_1,\dots,g_k}$ vanish, 
and hence $\omega\in L$. 

Next, we turn to the question of surjectivity. 
We shew that, when $p\leq n-1$ is a prime which is not inversible in $\B$, 
the elements $\zeta_1=e_1^{[p]}e_2e_3\cdots e_{n-p+1}$ and $\zeta_2= e_1e_2^{[p]}e_3\cdots e_{n-p+1}$ 
do not belong to $\Im\gamma_n$. Since the basis elements of $\BM$, by the lemma, 
transform under $\gamma_n$ as 
\beq				\label{E: Image of basis}
[e_{p_1}\de\cdots\de e_{p_m}] \mapsto 
\sum_{\substack{g_1+\cdots+g_m \\ g_i\geq 1}} e_{p_1}^{[g_1]}\cdots e_{p_m}^{[g_k]},
\eeq
the elements $\zeta_1$ and $\zeta_2$ can only be produced as images of elements of the form 
$$
\epsilon_{r,s} = \left[\De_{r} e_1\de \De_{s} e_2 \de e_3 \cdots \de e_{n-p+1}\right].
$$

If $s\geq 2$, then $\zeta_1$ does not occur in the expression for $\gamma_n(\epsilon_{r,s})$. 
When $s=1$, the part corresponding to $\zeta_1$ in 
$\gamma_n(\epsilon_{r,s})$ is, using \eref{E: Image of basis},  
$$
\sum_{\substack{g_1+\cdots+g_r=p \\ g_i\geq 1}} e_1^{[g_1]}\cdots e_1^{[g_r]} e_2 e_3\cdots e_{n-p+1} 
= \sum_{\substack{g_1+\cdots+g_r=p \\ g_i\geq 1}} 
\binom{p}{g_1,\dots, g_r} e_1^{[p]} e_2 e_3\cdots e_{n-p+1}.
$$
Each term is divisible by $p$ if $r\geq 2$.
(This may also be deduced from the fact that the multinomial sum 
counts surjective functions $[p]\to [r]$, 
and so the expression equals $r! \stir{p}{r} \zeta_1$.)
The situation is, of course, symmetric for $\zeta_2$.
Hence $\zeta_1$ and $\zeta_2$ only occur with non-zero co-efficients modulo $p$ in 
$\gamma_n(\epsilon_{1,1})$, in which they both occur with co-efficient $1$. Clearly, they cannot lie 
in $\Im\gamma_n$.

Conversely, suppose $(n-1)!$ is inversible. If $g_1+\cdots+g_k=n$, one has 
$$
\gamma_n \left(\left[ \De_{g_1} e_1 \de \cdots \de \De_{g_k} e_k \right] \right) 
= e_1^{g_1} \cdots e_k^{g_k} = g_1!\cdots g_k! e_1^{[g_1]} \cdots e_k^{[g_k]} \in \Im\gamma_n,
$$
where the co-efficient is inversible if all $g_i\leq n-1$. An exception occurs if only one $g_i=n$, 
the rest being $0$, but in this case 
$\gamma_n([\de e_i])=e_i^{[n]}\in \Im\gamma_n$.
\epr




\section{Numerical versus Strict Polynomial Functors}

We now set forth one or two comparison theorems for numerical and strict polynomial functors, 
sprung from properties of the natural functor $\HPol_n\to\Num_n$. These categories 
are, by Theorems \ref{S: HPol} and \ref{S: Num}, 
equivalent to $\GBMod$ and $\BBMod$, respectively.

To this end, we specialise the discussion in the preceding section to the particular case 
$M=\B^{n\times n}$. Both $\BB$ and $\GB$ are algebras, and 
$$
\gamma_n\colon \BB \to \GB 
$$ 
will be an homomorphism of such, since  
$$ 
\gamma_n([\alpha])\circ \gamma_n([\beta]) 
= \alpha^{[n]} \circ \beta^{[n]} = (\alpha\beta)^{[n]} 
= \gamma_n([\alpha\beta]) = \gamma_n([\alpha][\beta]). 
$$
Associated to this algebra homomorphism is the restriction of scalars 
$$ 
\GBMod\to\BBMod ,
$$
which views a $\GB$-module as a $\BB$-module under the action
$ [\sigma]x = \gamma_n(\sigma)x = \sigma^{[n]}x$.
On the level of polynomial functors, this corresponds to the forgetful functor $\HPol_n\to\Num_n$. 
As such, this functor will necessarily be exact and faithful 
(its left adjoint being extension of scalars): 

\bth 				\label{T: Forgetful} 
The forgetful functor $\HPol_n\to\Num_n$ corresponds to restriction of scalars 
$$
\GBMod\to\BBMod
$$ 
along the divided power map $\gamma_n \colon \BB \to \GB$.
This functor is exact and faithful. It is full and essentially injective when $n=0,1,2$.
\eth

\bpr
We make use of the fact that a necessary and sufficient condition for 
restriction of scalars to be a full functor is that the associated  
ring homomorphism be epimorphic. (See e.g. \cite{Popescu}, Proposition 16.3.)
According to Theorem \ref{S: Kernel}, $\gamma_n$ is surjective for $n=0,1,2$, and so
\emph{a fortiori} epimorphic. 

Suppose $\mu,\nu\colon \GB\to\End M$ are two $\GB$-module structures 
on a $\B$-module $M$, co-inciding as restricted $\BB$-modules. Then $\mu\gamma_n=\nu\gamma_n$, 
which implies $\mu=\nu$ by the surjectivity of $\gamma_n$. Consequently, restriction of scalars 
is full and essentially injective for $n=0,1,2$.
\epr

One may enquire for the exact conditions under which the restriction functor above is full. 
Since this is equivalent to asking when $\gamma_n\colon \BM\to\GM$ is an epimorphism of rings, 
the question is more delicate than that answered in Theorem \ref{S: Kernel} above.
We do not pursue this direction, but content ourselves with an account of 
what may happen in the case $n=3$:

\bex				\label{X: gamma_3}
We shew that $\gamma_3$ is not epimorphic, provided $2$ is not inversible in $\B$.
Let $\sigma_{ij}\colon \B^3\to\B^3$, for $i,j=1,2,3$, transform 
the $j$'th basis vector to the $i$'th basis vector, and all other basis vectors to $0$.  
These maps form a basis for~$\B^{3\times 3}$. 

The unity of $\Gamma^3(\B^{3\times 3})/2$ has a decomposition 
\begin{multline*}
1^{[3]} = (\sigma_{11}+\sigma_{22}+\sigma_{33})^{[3]} 
= \sum_i \sigma_{ii}^{[3]} 
+ \sum_{i\neq j} (\sigma_{ii}^{[2]}\sigma_{jj}+\sigma_{ii}\sigma_{ij}\sigma_{ji}) \\
+ \sum_{i\neq j} \sigma_{ii}\sigma_{ij}\sigma_{ji} 
+ (\sigma_{11}\sigma_{22}\sigma_{33} + \sigma_{12}\sigma_{23}\sigma_{31} 
+ \sigma_{13}\sigma_{32}\sigma_{21} )
+ (\sigma_{12}\sigma_{23}\sigma_{31} + \sigma_{13}\sigma_{32}\sigma_{21} )
\end{multline*}
into orthogonal idempotents.
Further factoring out the idempotents on the second row 
 produces a factor ring $\pi\colon\Gamma^3(\B^{3\times 3})\to R$,  
 which is spanned, over $\B/2$, by the images of all divided powers of the form
$$
\sigma_{ii}^{[3]}, \qquad \sigma_{ii}^{[2]}\sigma_{jj}, 
\qquad \sigma_{ii}^{[2]}\sigma_{ij}, \qquad \sigma_{ii}^{[2]}\sigma_{ji}, 
\qquad \sigma_{ij}^{[2]}\sigma_{ii}, \qquad \sigma_{ji}^{[2]}\sigma_{ii},
$$
and having unity (the image of)
$$
 \sum_i \sigma_{ii}^{[3]} + \sum_{i\neq j} \sigma_{ii}^{[2]}\sigma_{jj}.
$$
An involution $\phi$ of $R$ is provided by conjugation by 
$$
\omega = \sum_i \sigma_{ii}^{[3]} + \sum_{i\neq j} \sigma_{ij}^{[2]}\sigma_{ji}.
$$
One readily verifies the elements
\begin{align*}
\gamma_3([\de\sigma_{ij}]) &= \sigma_{ij}^{[3]} \\ 
\gamma_3([\sigma_{ij}\de\sigma_{kl}]) &= \sigma_{ij}^{[2]}\sigma_{kl} + \sigma_{ij}\sigma_{kl}^{[2]} \\
\gamma_3([\sigma_{ij}\de\sigma_{kl}\de\sigma_{mn}]) &= \sigma_{ij}\sigma_{kl}\sigma_{mn} 
\end{align*}
to be either $0$ in $R$ or invariant under the action of $\phi$. 
It follows that $\phi\pi\gamma_3=\pi\gamma_3$, so that $\gamma_3$ is not epimorphic. 

Consequently, restriction of scalars cannot be full in this case. 
Nor is it essentially injective, since we may define two $\Gamma^3(\B^{3\times 3})$-module structures on $R$, 
viz.~letting $\zeta\in\Gamma^3(\B^{3\times 3})$ act as left multiplication by $\zeta$ and by
$\phi(\zeta)=\omega\circ\zeta\circ\omega^{-1}$, respectively. 
These are indeed distinct, for 
$$
\sigma_{11}^{[2]} \sigma_{22} \circ \sigma_{11}^{[2]} \sigma_{22} = \sigma_{11}^{[2]} \sigma_{22},
\quad\text{but}\quad
\omega\circ \sigma_{11}^{[2]} \sigma_{22} \circ \omega^{-1} \circ \sigma_{11}^{[2]} \sigma_{22} = 
\sigma_{22}^{[2]} \sigma_{11} \circ \sigma_{11}^{[2]} \sigma_{22} = 0.
$$
The structures co-incide as restricted $P_3(\B^{3\times 3})$-modules, 
since $\phi\pi\gamma_3=\pi\gamma_3$ as per the above.
\eex

Recall from Theorem \ref{S: Num} that, under the Morita equivalence, 
the numerical functor $F$ corresponds to the $\BB$-module $F(\B^n)$, and  
that, if $F$ is strictly polynomial, this is in fact a $\GB$-module.

\bth 			\label{S: Quasi-homogeneous}
The numerical functor $F$, of degree $n$, is quasi-homogeneous of degree $n$ if and only if $F(\B^n)$
is a module over $\Im\gamma_n\cong \BB/\Ker\gamma_n$.  
\eth

\bpr 
The action of $\BB$ on $F(\B^n)$ is given by 
$[\sigma]x = F(\sigma)(x)$. 
The requirement that $\Ker\gamma_n$ annihilate $F(\B^n)$ is equivalent to demanding 
that $F$ itself vanish on 
$$
\Ker\gamma_n =\Q\otimes_{\Z} \gen{[r\sigma] - r^n[\sigma] | r\in\B, \ \sigma\in \B^{n\times n} } \cap \BB
$$ 
(Theorem \ref{S: Kernel}), which is quasi-homogeneity. 
\epr

Plainly, a quasi-homogeneous functor $F$ of degree $n$ 
admits an homogeneous structure if and only if the 
$\Im\gamma_n$-module structure on $F(\B^n)$ may be extended to 
a $\GB$-module structure. Quasi-homogeneity is necessary because 
$$
\BB/\Ker\gamma_n \cong \Im\gamma_n \subseteq \GB, 
$$
and so we have a necessary and sufficient condition for a numerical 
functor to be strictly polynomial:

\bth[The Polynomial Functor Theorem]			\label{T: PFT}
Let $F$ be a numerical functor of degree $n$. Then $F$ may be given the structure of homogeneous 
functor of degree~$n$ if and only if $F$ is quasi-homogeneous of degree $n$ and the 
 $\Im\gamma_n$-module structure on $F(\B^n)$ may be extended to 
a $\GB$-module structure. 
\eth

\bex 			\label{X: Quadratic} 
Any quasi-homogeneous functor of degree $2$  
may be given a unique strict polynomial structure, which makes it homogeneous of degree $2$. 
This stems from the map $\gamma_2$ being surjective, 
inducing an isomorphism 
$$
\B[\B^{2\times 2}]_2/\Ker\gamma_2\cong \Im\gamma_2=\Gamma^2(\B^{2\times 2}). \qedhere
$$
\eex

\section{The Hierarchy of Numerical Functors}

We proceed to discuss \emph{locally numerical} and \emph{analytic} functors. 
Salomonsson's investigations 
gave us the next result in the setting of strict polynomial functors (\cite{Pelle}, Propositions~2.3, 2.5), 
where he instead lets $A$ range over all \emph{commutative} algebras. 
In the case he considers, conditions A, B and C are equivalent if the arrow maps are assumed 
strictly polynomial of uniformly bounded degree. 

\bth 			\label{S: ABC}
Consider the following constructs, where $A$ ranges over all binomial algebras:
\balph
\item A family of ordinary functors $E_A\colon \AFMod\to \AMod$, commuting with extension of scalars. 
\item A functor $J\colon \FMod\to \Mod$ with arrow maps  
$$
J_A\colon \Hom_A(A\otimes M,A\otimes N) \to \Hom_A(A\otimes J(M),A\otimes J(N)),
$$ 
 multiplicative and natural in $A$. 
\item A functor $F\colon \FMod \to \Mod$ with arrow maps 
$$
F_A\colon A\otimes \Hom_\B(M,N) \to A\otimes \Hom_\B(F(M),F(N)),
$$
 multiplicative and natural in $A$ (a numerical functor as per Theorem \ref{S: BAlg-polynomial}). 
\ealph
Constructs A and B are equivalent, but weaker than C. If, in addition, the arrow maps are presumed 
numerical of (uniformly) bounded degree, all three are equivalent. 
\eth

\bpr
Given $E$, we may define $J$ by $J(M) = E_\B(M)$ and the diagram: 
$$ 
\xymatrix{
\Hom_A(A\otimes M, A\otimes N) \ar[r]^-{E_A} \dotar[dr]_-{J_A} & \Hom_A(E_A(A\otimes M),E_A(A\otimes N)) \isoar[d] \\
& \Hom_A(A\otimes E_\B(M),A\otimes E_\B(N)) 
} 
$$ 
Conversely, starting from $J$, the functors $E_A$ may be defined by the equation $E_A(M)=A\otimes J(M)$ and 
letting
$$ 
E_A = J_A \colon 
\xymatrix{
\Hom_A(A\otimes M,A\otimes N) \ar[r] & \Hom_A(A\otimes J(M),A\otimes J(N)).
} 
$$ 
Next, to define $J$ from $F$, let $J(M)=F(M)$ and use the diagram:
$$ 
\xymatrix{
A\otimes \Hom_\B(M,N) \isoar[d] \ar[r]^-{F_A} & A\otimes \Hom_\B(F(M),F(N)) \ar[d] \\
\Hom_A(A\otimes M,A\otimes N) \dotar[r]_-{J_A} & \Hom_A(A\otimes F(M),A\otimes F(N)) 
} 
$$ 
The left column is an isomorphism in so far as $M$ and $N$ are free. 

The difficult part is defining $F$ from $J$, provided that $J$ is indeed of bounded degree $n$. 
The proof is modelled on Salomonsson's argument for strict polynomial functors. 
We define $F=J$ on objects, and let
$$
F_\B = J_\B\colon \Hom_\B(M,N) \to \Hom_\B(J(M),J(N))
$$
for $\B$-modules $M$ and $N$. 
If this map can be shewn to be numerical, it may, by Theorem 10 of \cite{PM}, be uniquely extended to a 
family of maps $F_A$ as above, natural in the binomial algebra $A$.

To this end, find a free resolution
$$ 
\xymatrix{
\B^{(\lambda)} \ar[r] & \B^{(\kappa)} \ar[r] & J(M) \ar[r] & 0.
} 
$$
Applying the left-exact functor $\Hom_\B( -,J(N))$ yields a commutative diagram:
$$ 
\xymatrix{
0 \ar[r] & \Hom_\B(J(M),J(N)) \ar[r]^-\iota & J(N)^{\kappa} \ar[r]^-\sigma & J(N)^{\lambda} \\
&  \Hom_\B(M,N) \ar[u]^{J_\B} \ar[r]_-{\delta_n} & P_n(\Hom_\B(M,N)) \dotar[u]_\zeta \dotar[ul]_{\xi}
} 
$$
The homomorphism $\iota J_\B$ may be split up into components 
$$
(\iota J_\B)_k\colon \Hom_\B(M,N) \to J(N), \qquad k\in \kappa.
$$ 
These are, by the assumption on $J$, numerical of degree $n$, and will therefore factorise through $\delta_n$ 
via some linear $\zeta_k\colon P_n(\Hom_\B(M,N))\to J(N)$. This establishes the existence of a linear map $\zeta$ as in the diagram, 
making the square commute. 

Now, $\sigma\zeta\delta_n = \sigma\iota J_\B =0$, giving $\sigma\zeta=0$. By the exactness of the upper row, 
$\zeta$ admits a factorisation through some homomorphism $\xi$ as in the diagram. 
Because $ \iota J_\B = \zeta\delta_n = \iota\xi\delta_n$ 
and $\iota$ is one-to-one, $J_\B=\xi\delta_n$. 
Since $J_\B$ admits a factorisation through $P_n(\Hom_\B(M,N))$, it is indeed numerical of degree $n$.
\epr

Numerical functors satisfy all three conditions in the theorem. 
A functor will be called \textbf{analytic} if it only satisfies the weaker conditions A and B. 
Examples are the classical algebraic functors $T$ (the tensor algebra), $S$ 
(the symmetric algebra) and $\Gamma$ (the divided power algebra). 

A functor fulfilling condition C, but without any assumption on bounded degree, may be 
called \textbf{locally numerical}. 
An example would be $\Lambda$ (the exterior algebra). This is because, when $n>p$, the module 
$ \Lambda^n(\B^p) = 0$, and hence 
$$ 
\Lambda\colon \Hom(\B^p,\B^q) \to \Hom(\Lambda(\B^p),\Lambda(\B^q)) 
$$
is numerical of degree $\max(p,q)$.

\section{Analytic Functors}

A strict analytic functor is well known to be the inductive limit of 
its strict polynomial subfunctors. (In fact, it is a direct sum of such.)
There is no dissimilarity in the numerical case. 
Whereas, traditionally, analytic functors have been
identified with the inductive limits of polynomial functors, the definition we 
gave above is no different, as we now set out to shew.

\blem 
Let $F$ be an analytic functor and $P$ a finitely generated, free module. 
Suppose $u\in F(P)$, and define the subfunctor $G$ by 
$$ 
G(M) = \gen{ F(\alpha)(u) |  \alpha\colon P\to M}. 
$$
Consider the natural transformation $\xi\colon \Hom(P,-)\to F$, given by 
$$
\xi_N \colon \Hom(P,N) \to F(N), \qquad \alpha\mapsto F(\alpha)(u).
$$
If $\xi_N$, for some $N$, is numerical of degree $n$, 
then so is 
$$
G_{M,N}\colon \Hom(M,N)\to \Hom(G(M),G(N))
$$
for any $M$. In particular: 
\begin{itemize}
 \item If all $\xi_N$ are numerical, then $G$ is locally numerical. 
\item If all $\xi_N$ are numerical of uniformly bounded degree, then 
$G$ is numerical. 
\end{itemize}
\hfill \elem

\bpr 
The modules $G(M)$ are invariant under the action of $F$, so $G$ is indeed a subfunctor of $F$.
Suppose $\xi_N$ is numerical of degree $n$. Then, for all homomorphisms 
$\alpha,\alpha_i\colon P\to N$ and scalars $r$, the equations 
$$ 
F(\alpha_1\de \cdots \de \alpha_{n+1})(u) = 0 
\qquad\text{and}\qquad F(r\alpha)(u) = \sum_{m=0}^n \binom rm F\left(\De_m \alpha\right)(u) 
$$
hold. This implies that, for all scalars $r$ and homomorphisms  
$\beta,\beta_i\colon M\to N$ and  $\gamma\colon P\to M$, the following equations hold true:
$$ 
F(\beta_1\de \cdots \de \beta_{n+1})F(\gamma)(u) = 0 
\quad\text{and}\quad 
F(r\beta)F(\gamma)(u) = \sum_{m=0}^n \binom rm F\left(\De_m \beta\right)F(\gamma)(u). 
$$
Hence 
$$ 
F(\beta_1\de \cdots \de \beta_{n+1}) = 0 
\qquad\text{and}\qquad 
F(r\beta) = \sum_{m=0}^n \binom rm F\left(\De_m \beta\right)
$$
on $G(M)$, so that $G_{M,N}$ is indeed numerical of degree $n$. 
\epr

\bth 			\label{S: Analytic}
The analytic functors are precisely the filtered inductive limits of numerical functors. 
\eth

\bpr
\emph{Step 1: Filtered inductive limits of numerical, or even analytic, functors are analytic.}
Let $F_i$, for $i\in I$, be analytic functors, and let $A$ be a binomial algebra. Any 
$\alpha\in\Hom_A(A\otimes M,A\otimes N)$ gives rise to a map
$F_i(\alpha)\colon A\otimes F_i(M) \to A\otimes F_i(N)$.
Therefore 
$$ 
\indlim F_i(\alpha)\colon A\otimes \indlim F_i(M) \to A\otimes \indlim F_i(N)
$$
since tensor products commute with inductive limits, which yields a map 
$$ 
\indlim F_i \colon \Hom_A(A\otimes M, A\otimes N) 
\to \Hom_A(A\otimes \indlim F_i(M), A\otimes \indlim F_i(N)),  
$$
establishing that $\indlim F_i$ is analytic. 

\emph{Step 2: Analytic functors are filtered inductive limits of locally numerical functors.} 
Let $F$ be an analytic functor. The maps 
$$
F_{M,N}\colon \Hom_A(A\otimes M,A\otimes N) \to \Hom_A(A\otimes F(M),A\otimes F(N))
$$ 
are then multiplicative and natural in the binomial algebra $A$. 
To shew $F$ is the inductive limit of locally numerical functors,  
it is sufficient to construct, for any free and finitely generated 
module $P$ and any element $u\in F(P)$, 
a locally numerical subfunctor $G$ of $F$ such that $u\in G(P)$. 
To this end, define $G$ as in the lemma;
then clearly $u\in G(P)$. If only we can shew that 
$\xi_M\colon \Hom(P,M) \to F(M)$ is always numerical (of possibly unbounded degree), 
then $G$ will be locally numerical. 

To this end, let $\epsilon_1,\dots,\epsilon_k$ be a basis of 
the free and finitely generated module $\Hom(P,M)$.
Let
$ A= \freenum[\B]{s_1,\dots,s_k}$ be the free binomial algebra on the variables $s_1,\dots,s_k$ 
(see \cite{BR}).
Since 
$$ 
F\left(\sum s_i \otimes \epsilon_i\right) \in \Hom_A(A\otimes F(P),A\otimes F(M)),
$$
we may write 
$$ 
F\left(\sum s_i \otimes \epsilon_i\right)(1 \otimes u) 
= \sum_{\substack{g_1+\cdots+g_k\leq n \\ g_i\geq 0}} 
\binom{s_1}{g_1}\cdots\binom{s_k}{g_k} \otimes v_{g_1,\dots,g_k} \in A\otimes F(M) 
$$
for some $n$. By the naturality of $F$, we may specialise each $s_i$ to an $a_i\in\B$:
$$ 
\xi_M\left(\sum a_i \epsilon_i\right) =  F\left(\sum a_i \epsilon_i\right)(u) 
= \sum_{g_1+\cdots+g_k\leq n} \binom{a_1}{g_1}\cdots\binom{a_k}{g_k} v_{g_1,\dots,g_k}.  
$$
Using Theorem 8 of \cite{PM}, it follows that $\xi_M$ is numerical of degree $n$. 

\emph{Step 3: Locally numerical functors are filtered inductive limits of numerical functors.}
Let $F$ be a locally numerical functor. Once again, given a free and finitely generated module $P$ 
and $u\in F(P)$, define $G$ and $\xi$ as in the lemma. 
We shall shew that $G$ is numerical by shewing that all $\xi_M$ are numerical of 
uniformly bounded degree.  

Consider the free binomial algebras
$B = \freenum[\B]{s_1,\dots,s_k}$ and $C = \freenum[\B]{s_1,\dots,s_k,t}$.
Let $\tau\colon B\to C$ be the algebra homomorphism given by $s_i \mapsto ts_i$.
By the local numericality of $F$, there is a  commutative diagram: 
$$ 
\xymatrix{
B \ar[d]_{\tau} & B\otimes \Hom(P,M) \ar[d]_{\tau\otimes 1} \ar[r]^-F & B\otimes \Hom(F(P), F(M)) \ar[d]^{\tau\otimes 1} \\
C & C\otimes \Hom(P,M) \ar[r]_-F & C\otimes \Hom( F(P), F(M))
} 
$$
As a consequence, we obtain, for any homomorphisms $\alpha_i\colon P\to M$:
\beq				\label{E: Analytic Eq 1}
(\tau\otimes 1) F\left(\sum s_i \otimes \alpha_i\right) 
=  F\left((\tau\otimes 1)\left(\sum s_i \otimes \alpha_i\right)\right) 
= F\left(\sum ts_i\otimes\alpha_i \right).
\eeq

Considering now 
$$ 
F\colon \freenum[\B]{s_1,\dots,s_k}\otimes\Hom(P,M) 
\to \freenum[\B]{s_1,\dots,s_k}\otimes \Hom(F(P),F(M)),
$$
we write 
\beq				\label{E: Analytic Eq 2} 
F\left(\sum s_i \otimes \alpha_i\right) = \sum_{g_1,\dots,g_k\geq 0} 
\binom{s_1}{g_1}\cdots\binom{s_k}{g_k}  \otimes \beta_{g_1,\dots,g_k},
\eeq
for homomorphisms $\beta_{g_1,\dots,g_k}\colon F(P)\to F(M)$. Only a finite number of
these are non-zero, but this number depends on the $\alpha_i$.
Similarly, from contemplating
$$ 
F\colon \freenum[\B]{t}\otimes\Hom(P,P) \to \freenum[\B]{t}\otimes \Hom(F(P),F(P)),
$$
we may write 
\beq				\label{E: Analytic Eq 3} 
F(t\otimes 1_P) = \sum_{m\leq n} \binom tm \otimes \gamma_m,
\eeq
for some number $n$, depending purely on $P$, and homomorphisms $\gamma_m\colon F(P)\to F(P)$. 
 Now combine \eref{E: Analytic Eq 1}, \eref{E: Analytic Eq 2} 
and \eref{E: Analytic Eq 3}:
\begin{multline*}
\sum_{g_1,\dots,g_k} \binom {ts_1}{g_1}\cdots \binom{ts_k}{g_k} \otimes \beta_{g_1,\dots,g_k}
= (\tau\otimes 1) \left(\sum_{g_1,\dots,g_k} \binom{s_1}{g_1}\cdots\binom{s_k}{g_k} 
\otimes \beta_{g_1,\dots,g_k} \right)  \\
= (\tau\otimes 1) F\left(\sum s_i \otimes \alpha_i\right) 
= F\left( \sum ts_i \otimes \alpha_i \right) 
= F\left( \sum s_i \otimes \alpha_i \right)F(t\otimes 1_P) \\
= \left( \sum_{g_1,\dots,g_k} \binom{s_1}{g_1}\cdots\binom{s_k}{g_k} 
\otimes \beta_{g_1,\dots,g_k} \right) 
\left( \sum_{m\leq n} \binom tm \otimes \gamma_m \right) \\
= \sum_{g_1,\dots,g_k} \sum_{m\leq n} \binom{s_1}{g_1}\cdots\binom{s_k}{g_k} \binom tm
\otimes \beta_{g_1,\dots,g_k}\gamma_m.
\end{multline*}
The right-hand side, and therefore also the left-hand side, is of degree $n$ in $t$, whence 
$\beta_{g_1,\dots,g_k}=0$ when the sum of the indices exceeds $n$. Consequently, 
$$ 
\xi_M\left(\sum a_i \alpha_i\right) 
= F\left(\sum a_i \alpha_i\right)(u) 
= \sum_{g_1+\cdots+g_k\leq n} \binom{a_1}{g_1}\cdots\binom{a_k}{g_k} \beta_{g_1,\dots,g_k}(u), 
$$
where we have specialised $s_i$ to $a_i\in\B$, and so $\xi_M$ is numerical of degree $n$.
\epr

\end{document}